\documentclass[11pt]{article} 
\usepackage{amsmath,amssymb,amsthm}  

\allowdisplaybreaks

\begin{document}

\def\fl#1{\left\lfloor#1\right\rfloor}
\def\cl#1{\left\lceil#1\right\rceil}
\def\ang#1{\left\langle#1\right\rangle}
\def\stf#1#2{\left[#1\atop#2\right]} 
\def\sts#1#2{\left\{#1\atop#2\right\}}
\def\N{\mathbb N}
\def\Z{\mathbb Z}
\def\R{\mathbb R}
\def\C{\mathbb C}

\newtheorem{theorem}{Theorem}
\newtheorem{Prop}{Proposition}
\newtheorem{Cor}{Corollary}
\newtheorem{Lem}{Lemma}

\newenvironment{Rem}{\begin{trivlist} \item[\hskip \labelsep{\it
Remark.}]\setlength{\parindent}{0pt}}{\end{trivlist}}

\title{Weighted Sylvester sums on the Frobenius set
}

\author{
Takao Komatsu\, and Yuan Zhang\\
\small Department of Mathematical Sciences, School of Science\\[-0.8ex]
\small Zhejiang Sci-Tech University\\[-0.8ex]
\small Hangzhou 310018 China\\[-0.8ex]
\small \texttt{komatsu@zstu.edu.cn}\, and \texttt{zjlgzy51@126.com}
}

\date{
\small MR Subject Classifications: Primary 05A15; Secondary 11D07, 11B68, 11P81, 05A17 
}

\maketitle
 
\begin{abstract} 
Let $a$ and $b$ be relatively prime positive integers. 
In this paper the weighted sum $\sum_{n\in{\rm NR}(a,b)}\lambda^{n-1}n^m$ is given explicitly or in terms of the Apostol-Bernoulli numbers, where $m$ is a nonnegative integer, and ${\rm NR}(a,b)$ denotes the set of positive integers nonrepresentable in terms of $a$ and $b$. 
\\
{\bf Keywords:} Frobenius problem, weighted sums, Sylvester sums, Apostol-Bernoulli numbers      
\end{abstract}

\section{Introduction}  

The {\it Frobenius Problem} is to determine the largest positive integer that is NOT representable as a nonnegative integer combination of given positive integers that are coprime (see \cite{ra05} for general references). 

Given positive integers $a_1,\dots,a_m$ with $\gcd(a_1,\dots,a_m)=1$, it is well-known that all sufficiently large $n$ the equation 
\begin{equation}
a_1 x_1+\cdots+a_m x_m=n 
\label{fb:eq} 
\end{equation} 
has a solution with nonnegative integers $x_1,\dots,x_m$.   

The {\it Frobenius number} $F(a_1,\dots,a_m)$ is the LARGEST integer $n$ such that (\ref{fb:eq}) has no solution in nonnegative integers.  
For $m=2$, we have 
$$ 
F(a,b)=(a-1)(b-1)-1
$$ 
(Sylvester (1884) \cite{sy1884}).  For $m\ge 3$, exact determination of the Frobenius number is difficult.  
The Frobenius number cannot be given by closed formulas of a certain type (Curtis (1990) \cite{cu90}), the problem to determine $F(a_1,\dots,a_m)$ is NP-hard under Turing reduction (see, e.g., Ram\'irez Alfons\'in \cite{ra05}). Nevertheless, the Frobenius number for some special cases are calculated (e.g., \cite{op08,ro56,se77}). One convenient formula is by Johnson \cite{jo60}. One analytic approach to the Frobenius number can be seen in \cite{bgk01,ko03}.    
Some formulae for the Frobenius number in three variables can be seen in \cite{tr17}. 

For given $a$ and $b$ with $\gcd(a,b)=1$, let ${\rm NR}(a,b)$ denote the set of nonnegative integers nonrepresentable in term of $a$ and $b$, namely the set of all those nonnegative integers $n$ which cannot be expressed in the form $n=a x+b y$, where $x$ and $y$ are nonnegative integers. 

There are many kinds of problems related to the Frobenius problem. The problems for the number of solutions (e.g., \cite{tr00}), and the sum of integer powers of the gaps values in numerical semigroups (e.g., \cite{bs93,fr07,fks}) are popular. 
One of other famous problems is about the so-called {\it Sylvester sums} $\sum_{n\in{\rm NR}(a,b)}n^m$, where $m$ is a nonnegative integer (see, e.g., \cite{tu06} and references therein).   
Recently in \cite{bb20}, one of more general cases is considered by giving the largest integer, the number of integers and the sum of integers whose number of representation is less than or equal to $k$. When $k=1$, the problem is reduced to the classical case.

In this paper, we consider the weighted sum  
$$
S_m^{(\lambda)}(a,b):=\sum_{n\in{\rm NR}(a,b)}\lambda^{n-1}n^m\quad(\lambda\ne 0)\,. 
$$ 
Sylvester \cite{sy1884} showed that $S_0^{(1)}(a,b)=(a-1)(b-1)/2$, and 
Brown and Shuie showed \cite{bs93} that 
$$
S_1^{(1)}(a,b)=\frac{1}{12}(a-1)(b-1)(2 a b-a-b-1)\,. 
$$  
R\o dseth \cite{ro94} obtained a general formula for $S_m^{(1)}$ in terms of Bernoulli numbers and deduced 
$$
S_2^{(1)}(a,b)=\frac{1}{12}(a-1)(b-1)a b(a b-a-b)\,. 
$$ 
Tuenter \cite{tu06} also investigated $S_m^{(1)}$ by taking a different approach. He established relations between Sylvester sums and the power sums over the natural numbers. 
Wang and Wang \cite{wang08} considered the alternating Sylvester sums 
$$
T_m(a,b)=\sum_{n\in{\rm NR}(a,b)}(-1)^n n^m
$$ 
by using Bernoulli and Euler numbers.  

The purpose of this paper is to give an explicit expression for $S_m^{(\lambda)}(a,b)$. For $m=1$, we can give the following formula.   

\begin{theorem}  
For $\lambda\ne 0$ with $\lambda^a\ne 1$ and $\lambda^b\ne 1$, 
$$
S_1^{(\lambda)}(a,b)=\frac{1}{(\lambda-1)^2}+\frac{a b \lambda^{a b-1}}{(\lambda^a-1)(\lambda^b-1)}-\frac{(\lambda^{a b}-1)\bigl((a+b)\lambda^{a+b}-a \lambda^a-b \lambda^b\bigr)}{\lambda(\lambda^a-1)^2(\lambda^b-1)^2}\,. 
$$ 
\label{th:1}
\end{theorem}

We also give a general expression of $S_m^{(\lambda)}(a,b)$ in terms of the Apostol-Bernoulli numbers. The alternating Sylvester sums in \cite{wang08} can be also expressed as $T_m(a,b)=-S_m^{(-1)}(a,b)$.

The main new results (Theorems \ref{th:m} and \ref{th:m1} below) cover all values of $m$ and $\lambda$, and express $S_m^{(\lambda}(a,b)$ in terms of the Apostol-Bernoulli numbers. 
In case $m=1$ and $\lambda^a\not=1$ the expressions reduce to those given explicitly in Theorem \ref{th:1}.

\section{An explicit expression for $m=1$}  

As in \cite{bs93}, define 
$$
f(x)=\sum_{n=0}^{a b-a-b}\bigl(1-r(n)\bigr)x^n\,,
$$ 
where $r(n)$ denotes the number of representations of $n$ in the form $n=s a+t b$, where $s$ and $t$ are nonnegative integers.  
Since $r(n)=0$ or $1$ for $0\le n\le a b-1$, we have 
\begin{align*}
f'(\lambda)&=\sum_{n=1}^{a b-a-b}n\bigl(1-r(n)\bigr)\lambda^{n-1}=\sum_{1\le n\le a b -a-b\atop r(n)=0}n\lambda^{n-1}\\
&=\sum_{n\in{\rm NR}(a,b)}\lambda^{n-1}n=S_1^{(\lambda)}(a,b)\,.
\end{align*}  
We use the following fact from \cite{bs93}.  

\begin{Lem}
$$
f(x)=\frac{g(x)}{h(x)}\,, 
$$ 
where 
$$
g(x)=\sum_{k=1}^{b-1}\frac{x^{a k}-x^k}{x-1}\quad\hbox{and}\quad h(x)=\sum_{k=0}^{b-1}x^k\,.
$$ 
\label{lem-bs}
\end{Lem}

Suppose that $\lambda\ne 1\ne \lambda^a$. Then 
$$
h(\lambda)=\frac{\lambda^b-1}{\lambda-1}
$$ 
and 
\begin{align*}
h'(\lambda)=\sum_{k=0}^{b-1}k \lambda^{k-1}
&=\frac{b \lambda^{b-1}}{\lambda-1}-\frac{\lambda^b-1}{(\lambda-1)^2}\,. 
\end{align*}
Also, we have 
$$
g(\lambda)=\frac{(\lambda^{a b}-1)(\lambda-1)-(\lambda^a-1)(\lambda^b-1)}{(\lambda^a-1)(\lambda-1)^2}
$$ 
and 
\begin{align*}
g'(\lambda)&=\frac{(a b+1)\lambda^{a b}-a b \lambda^{a b-1}-(a+b)\lambda^{a+b+1}+a \lambda^{a-1}+b \lambda^{b-1}-1}{(\lambda^a-1)(\lambda-1)^2}\\
&\quad -\frac{a \lambda^{a-1}}{\lambda^a-1}g(\lambda)-\frac{2}{\lambda-1}g(\lambda)\,. 
\end{align*}
Hence, we finally get 
\begin{align*}
S_1^{(\lambda)}(a,b)&=f'(\lambda)=\frac{g'(\lambda)h(\lambda)-g(\lambda)h'(\lambda)}{\bigl(h(\lambda)\bigr)^2}\\ 
&=\frac{1}{(\lambda-1)^2}+\frac{a b \lambda^{a b-1}}{(\lambda^a-1)(\lambda^b-1)}
-\frac{(\lambda^{a b}-1)\bigl((a+b)\lambda^{a+b}-a \lambda^a-b \lambda^b\bigr)}{\lambda(\lambda^a-1)^2(\lambda^b-1)^2}\,. 
\end{align*}

In particular, for $\lambda=2$, we have the following. 

\begin{Cor}
\begin{align*}
\sum_{n\in{\rm NR}(a,b)}2^{n-1}n&=1+\frac{a b 2^{a b-1}}{(2^a-1)(2^b-1)}\\
&\quad -\frac{(2^{a b}-1)\bigl((a+b)2^{a+b}-2^a a-2^b b\bigr)}{2(2^a-1)^2(2^b-1)^2}\,. 
\end{align*}
\end{Cor} 

For example, for $a=3$ and $b=17$, 
\begin{align*} 
S_1^{(2)}(3,17)&=2^0\cdot 1+2^1\cdot 2+2^3\cdot 4+2^4\cdot 5+2^6\cdot 7+2^7\cdot 8+2^9\cdot 10\\
&\quad +2^{10}\cdot 11+2^{12}\cdot 13+2^{13}\cdot 14+2^{15}\cdot 16+2^{18}\cdot 19+2^{21}\cdot 22\\
&\quad +2^{24}\cdot 25+2^{27}\cdot 28+2^{30}\cdot 31\\
&=37515351605\,. 
\end{align*}
From Theorem \ref{th:1} (or the above Corollary), 
\begin{align*} 
S_1^{(2)}(3,17)&=\frac{1}{(2-1)^2}+\frac{3\cdot 17\cdot 2^{3\cdot 17-1}}{(2^3-1)(2^{17}-1)}\\
&\quad -\frac{(2^{3\cdot 17}-1)\bigl((3+17)2^{3+17}-3\cdot 2^3-17\cdot 2^{17}\bigr)}{2(2^3-1)^2(2^{17}-1)^2}\\
&=37515351605\,. 
\end{align*}

Similarly, by replacing $2$ by another value, we can obtain that  
\begin{align*} 
S_1^{(5)}(3,17)&=900879734470832437423896\,,\\
S_1^{(1/2)}(3,17)&=\frac{8822132865}{1073741824}\,,\\
S_1^{(-1)}(3,17)&=408\,,\\
S_1^{(-5/3)}(3,17)&=\frac{760508529478902941119864}{205891132094649}\,,\\
S_1^{(\pm\sqrt{2})}(3,17)&=34250061\pm 6965604\sqrt{2}\,. 
\end{align*}

\section{Weighted sums of higher power}  

Since 
\begin{align*}
f''(x)&=\frac{g''(x)}{h(x)}-\frac{2 g'(x)h'(x)+h(x)''(x)}{\bigl(h(x)\bigr)^2}+\frac{2 g(x)\bigl(h'(x)\bigr)^2}{\bigl(h(x)\bigr)^3}\\
&=\sum_{n=2}^{a b-a-b}n(n-1)\bigl(1-r(n)\bigr)x^{n-2}\,, 
\end{align*}
we get 
$$
x f''(x)+f'(x)=\sum_{n=0}^{a b-a-b}n^2\bigl(1-r(n)\bigr)x^{n-1}\,. 
$$ 
Hence, 
$$
S_2^{(\lambda)}(a,b)=\lambda f''(\lambda)+f'(\lambda)\,.
$$ 

For simplicity, put $X_1=(a+b)\lambda^{a+b}-a\lambda^a-b\lambda^b$ and $X_2=(a+b)^2\lambda^{a+b}-a^2\lambda^a-b^2\lambda^b$. Since 
$$
f'(\lambda)=\frac{1}{(\lambda-1)^2}+\frac{a b \lambda^{a b-1}}{(\lambda^a-1)(\lambda^b-1)}-\frac{(\lambda^{a b}-1)X_1}{\lambda(\lambda^a-1)^2(\lambda^b-1)^2}\,, 
$$ 
we get 
\begin{align*}
f''(\lambda)&=-\frac{2}{(\lambda-1)^3}+\frac{a b(a b-1)\lambda^{a b-2}}{(\lambda^a-1)(\lambda^b-1)}-\frac{2 a b\lambda^{a b-2}X_1}{(\lambda^a-1)^2(\lambda^b-1)^2}\\
&\quad-\frac{(\lambda^{a b}-1)(X_2-X_1)}{\lambda^2(\lambda^a-1)^2(\lambda^b-1)^2}+\frac{2(\lambda^{a b}-1)X_1}{\lambda^3(\lambda^a-1)^3(\lambda^b-1)^3}\,.
\end{align*}
Therefore, we obtain 
\begin{align*} 
S_2^{(\lambda)}(a,b)&=-\frac{\lambda+1}{(\lambda-1)^2}+\frac{a^2 b^2\lambda^{a b-1}}{(\lambda^a-1)(\lambda^b-1)}-\frac{2 a b\lambda^{a b}X_1+(\lambda^{a b}-1)X_2}{\lambda(\lambda^a-1)^2(\lambda^b-1)^2}\\
&\quad+\frac{2(\lambda^{a b}-1)X_1}{\lambda^2(\lambda^a-1)^3(\lambda^b-1)^3}\,.
\end{align*}

Similarly, we see that 
\begin{align*}
S_3^{(\lambda)}(a,b)&=\lambda^2 f'''(\lambda)+3 \lambda f''(\lambda)+f'(\lambda)\,,\\ 
S_4^{(\lambda)}(a,b)&=\lambda^3 f^{(4)}(\lambda)+6\lambda^2 f'''(\lambda)+7 \lambda f''(\lambda)+f'(\lambda)\,,\\ 
S_5^{(\lambda)}(a,b)&=\lambda^4 f^{(5)}(\lambda)+10 \lambda^3 f^{(4)}(\lambda)+25 \lambda^2 f'''(\lambda)+15\lambda f''(\lambda)+f'(\lambda)\,. 
\end{align*}

\section{Apostol-Bernoulli numbers}  

Though we may obtain explicit expressions of $S_m^{(\lambda)}(a,b)$ for small positive integers $m$, it is hard to obtain the formula for large $m$. In this section, using the so-called Apostol-Bernoulli numbers, we give an expression of $S_m^{(\lambda)}(a,b)$ for general positive integral $m$.  

Apostol-Bernoulli polynomials $\mathcal B_n(x,\lambda)$ are defined by the generating function (\cite[p.165, (3.1)]{ap51}) 
\begin{equation}  
\frac{z e^{x z}}{\lambda e^z-1}=\sum_{n=0}^\infty\mathcal B_n(x,\lambda)\frac{z^n}{n!}\quad(|z+\log\lambda|<2\pi)\,.  
\label{def:aposberpoly} 
\end{equation} 
When $\lambda=1$ in (\ref{def:aposberpoly}), $B_n(x)=\mathcal B_n(x,1)$ are the classical Bernoulli numbers.  When $x=0$ in (\ref{def:aposberpoly}), $\mathcal B_n(\lambda)=\mathcal B_n(0,\lambda)$ are {\it Apostol-Bernoulli numbers} (\cite[Definition 1.2]{luo04}), defined by 
\begin{equation}  
\frac{z}{\lambda e^z-1}=\sum_{n=0}^\infty\mathcal B_n(\lambda)\frac{z^n}{n!}\quad(|z+\log\lambda|<2\pi)\,.  
\label{def:aposbernum} 
\end{equation} 
They seem to be also called $\lambda$-Bernoulli numbers. 
When $\lambda=1$, the generating function of the left-hand side in (\ref{def:aposbernum}) is exactly the same as that of the classical Bernoulli numbers $B_n$. But it does not imply that $\mathcal B_n(1)=B_n$ on the right-hand side though quite a few authors misunderstand. In fact, as seen in \cite[p.165]{ap51}, the first several values are given by 
\begin{align*} 
&\mathcal B_0(\lambda)=0,\quad \mathcal B_1(\lambda)=\frac{1}{\lambda-1},\quad \mathcal B_2(\lambda)=-\frac{2\lambda}{(\lambda-1)^2},\quad \mathcal B_3(\lambda)=\frac{3\lambda(\lambda+1)}{(\lambda-1)^3},\\
&\mathcal B_4(\lambda)=-\frac{4\lambda(\lambda^2+4\lambda+1)}{(\lambda-1)^4},\quad \mathcal B_5(\lambda)=\frac{5\lambda(\lambda^3+11\lambda^2+11\lambda+1)}{(\lambda-1)^5}\,. 
\end{align*}
But,  
$$
B_0=1,~B_1=-\frac{1}{2},~B_2=\frac{1}{6},~B_3=0,~B_4=-\frac{1}{30},~B_5=0,~B_6=\frac{1}{42},~\dots\,. 
$$     

For $\lambda\ne 1$, Apostol Bernoulli polynomials $\mathcal B_n(x,\lambda)$ can be expressed explicitly by 
\begin{equation}  
\mathcal B_n(x,\lambda)=\sum_{k=1}k\binom{n}{k}\sum_{j=0}^{k-1}(-1)^j\lambda^j(\lambda-1)^{-j-1}j!\sts{k-1}{j}x^{n-k}\quad(n\ge 0)  
\label{exp:aposberpoly}
\end{equation}
(\cite[Remark 2.6]{luo04}), where the Stirling numbers of the second kind $\sts{n}{k}$ are given by 
$$
\sts{n}{k}=\frac{1}{k!}\sum_{j=0}^k(-1)^{k-j}\binom{k}{j}j^n\,. 
$$    
When $x=0$ in (\ref{exp:aposberpoly}), Apostol-Bernoulli numbers $\mathcal B_n(\lambda)$ have an explicit expression in terms of the Stirling numbers of the second kind (\cite[p.166, (3.7)]{ap51},\cite[p.510, (3)]{luo04})\footnote{In both references, the sum begins from $j=1$. However, the value for $n=1$ does not match the correct one $\mathcal B_1(\lambda)=1/(\lambda-1)$.}.  
\begin{equation}  
\mathcal B_n(\lambda)=n\sum_{j=0}^{n-1}(-1)^j\lambda^j(\lambda-1)^{-j-1}j!\sts{n-1}{j}\quad(n\ge 0)  
\label{exp:aposbernum}
\end{equation}

We use the similar approach by R\o dseth in (\cite{ro94}). 
Let $n$, $r$ and $s$ be integers with 
$$
r\equiv n\pmod a\quad (0\le r<a),\qquad b s\equiv r\pmod a \quad (0\le s<a)\,. 
$$ 
Notice that 
\begin{align*}
n\in{\rm NR}(a,b)&\Longleftrightarrow \exists t\in\mathbb Z~(1\le t\le\fl{b s/a}),~ n=-at +b s\\
&\Longleftrightarrow \exists k\in\mathbb Z~(0\le k\le(b s-r)/a-1),~ n=a k+r\,. 
\end{align*} 
Note that the case $\lambda=1$ is discussed in (\cite{ro94}).  
Since 
$$
S_m^{(\lambda)}(a,b)=\sum_{r=0}^{a-1}\sum_{k=0}^{\frac{b s-r}{a}-1}\lambda^{a k+r-1}(a k+r)^m\,, 
$$
for $\lambda\ne 1$, we have 
\begin{align}  
\sum_{m=0}^\infty S_m^{(\lambda)}(a,b)\frac{z^m}{m!}&=\frac{1}{\lambda}\sum_{r=0}^{a-1}\sum_{k=0}^{\frac{b s-r}{a}-1}(\lambda e^z)^{a k+r}\notag\\
&=\frac{1}{\lambda}\frac{1}{(\lambda e^z)^a-1}\left(\sum_{r=0}^{a-1}(\lambda e^z)^{b s}-\sum_{r=0}^{a-1}(\lambda e^z)^r\right)\notag\\
&=\frac{1}{\lambda}\frac{1}{(\lambda e^z)^a-1}\left(\sum_{s=0}^{a-1}(\lambda e^z)^{b s}-\sum_{r=0}^{a-1}(\lambda e^z)^r\right)\notag\\
&=\frac{1}{\lambda}\frac{a z}{(\lambda e^z)^a-1}\frac{b z}{(\lambda e^z)^b-1}\frac{(\lambda e^z)^{a b}-1}{a b z^2}-\frac{1}{\lambda}\frac{1}{\lambda e^z-1}\,.
\label{eq:gf}
\end{align}
Assume that $\lambda^a\ne 1$ and $\lambda^b\ne 1$. 
The second term (without sign) of the right-hand side is equal to 
\begin{align*}  
\frac{1}{\lambda}\frac{1}{\lambda e^z-1}&=\frac{1}{\lambda z}\sum_{m=0}^\infty\mathcal B_m(\lambda)\frac{z^m}{m!}\\
&=\frac{1}{\lambda}\sum_{m=0}^\infty\frac{\mathcal B_m(\lambda)}{m}\frac{z^{m-1}}{(m-1)!}\\
&=\frac{1}{\lambda}\sum_{m=0}^\infty\frac{\mathcal B_{m+1}(\lambda)}{m+1}\frac{z^{m}}{m!}\quad (\mathcal B_0(\lambda)=0)\,. 
\end{align*}
The first term is divided into two parts. One part (without sign) is given as 
\begin{align*}  
&\frac{1}{\lambda}\frac{1}{a b z^2}\frac{a z}{(\lambda e^z)^a-1}\frac{b z}{(\lambda e^z)^b-1}\\
&=\frac{1}{\lambda}\frac{1}{a b z^2}\left(\sum_{i=0}^\infty\mathcal B_i(\lambda^a)a^i\frac{z^i}{i!}\right)\left(\sum_{j=0}^\infty\mathcal B_j(\lambda^b)b^i\frac{z^j}{j!}\right)\\
&=\frac{1}{\lambda}\sum_{m=0}^\infty\sum_{i=0}^m\binom{m}{i}a^{i-1}b^{m-i-1}\mathcal B_i(\lambda^a)\mathcal B_{m-i}(\lambda^b)\frac{z^{m-2}}{m!}\\ 
&=\frac{1}{\lambda}\sum_{m=0}^\infty\frac{1}{(m+1)(m+2)}\sum_{i=0}^{m+2}\binom{m+2}{i}a^{i-1}b^{m-i+1}\mathcal B_i(\lambda^a)\mathcal B_{m-i+2}(\lambda^b)\frac{z^{m}}{m!}\,. 
\end{align*}
Another part is given as 
\begin{align*}  
&\frac{\lambda^{a b-1}}{a b z^2}\frac{a z}{(\lambda e^z)^a-1}\frac{b z}{(\lambda e^z)^b-1}e^{a b z}\\
&=\lambda^{a b-1}\left(\sum_{k=0}^\infty a^k b^k\frac{z^k}{k!}\right)\\
&\quad\times\left(\sum_{\ell=0}^\infty\frac{1}{(\ell+1)(\ell+2)}\sum_{i=0}^{\ell+2}\binom{\ell+2}{i}a^{i-1}b^{\ell-i+1}\mathcal B_i(\lambda^a)\mathcal B_{\ell-i+2}(\lambda^b)\frac{z^{\ell}}{\ell!}\right)\\
&=\lambda^{a b-1}\sum_{m=0}^\infty\sum_{\ell=0}^m\binom{m}{\ell}\frac{1}{(\ell+1)(\ell+2)}\\
&\quad\times\sum_{i=0}^{\ell+2}\binom{\ell+2}{i}a^{m-\ell+i-1}b^{m-i+1}\mathcal B_i(\lambda^a)\mathcal B_{\ell-i+2}(\lambda^b)\frac{z^m}{m!}\,.  
\end{align*}
Comparing the coefficients on both sides of (\ref{eq:gf}), we get the following expression.  

\begin{theorem}  
For $\lambda\ne 0$ with $\lambda^a\ne 1$ and $\lambda^b\ne 1$, and a nonnegative integer $m$, 
\begin{align*}
S_m^{(\lambda)}(a,b)&=\lambda^{a b-1}\sum_{\ell=0}^m\sum_{i=0}^{\ell+2}\binom{\ell+2}{i}\binom{m}{\ell}\frac{a^{m-\ell+i-1}b^{m-i+1}}{(\ell+1)(\ell+2)}\mathcal B_i(\lambda^a)\mathcal B_{\ell-i+2}(\lambda^b)\\
&\quad -\frac{1}{(m+1)(m+2)\lambda}\sum_{i=0}^{m+2}\binom{m+2}{i}a^{i-1}b^{m-i+1}\mathcal B_i(\lambda^a)\mathcal B_{m-i+2}(\lambda^b)\\ 
&\quad -\frac{\mathcal B_{m+1}(\lambda)}{(m+1)\lambda}\,. 
\end{align*}
\label{th:m}
\end{theorem}

\noindent 
{\it Remark.}  
When $m=1$ in the expression of Theorem \ref{th:m}, that of Theorem \ref{th:1} is reduced.
\medskip

If $\lambda^a=1$ or $\lambda^b=1$ in (\ref{eq:gf}), without loss of generality, we can assume that $\lambda^a=1$ and $\lambda^b\ne 1$. Because $\gcd(a,b)=1$,  $\lambda^a=\lambda^b=1$ is impossible for $\lambda\ne 1$. 
Then, the first term of the right-hand side of (\ref{eq:gf}) is equal to 
\begin{align*}
&\frac{1}{\lambda}\frac{a z}{e^{a z}-1}\frac{b z}{\lambda^b e^{b z}-1}\frac{e^{a b z}-1}{a b z^2}\\
&=\frac{1}{\lambda z}\left(\sum_{k=0}^\infty\frac{a^k b^k}{k+1}\frac{z^k}{k!}\right)\left(\sum_{i=0}^\infty B_i a^i\frac{z^i}{i!}\right)\left(\sum_{j=0}^\infty\mathcal B_j(\lambda^b) b^j\frac{z^j}{j!}\right)\\
&=\frac{1}{\lambda z}\left(\sum_{k=0}^\infty\frac{a^k b^k}{k+1}\frac{z^k}{k!}\right)\left(\sum_{\ell=0}^\infty\sum_{i=0}^\ell\binom{\ell}{i}a^i b^{\ell-i}B_i\mathcal B_{\ell-i}(\lambda^b)\frac{z^\ell}{\ell!}\right)\\
&=\frac{1}{\lambda z}\sum_{m=0}^\infty\sum_{\ell=0}^m\sum_{i=0}^\ell\binom{m}{\ell}\binom{\ell}{i}\frac{a^{m-l+i}b^{m-i}}{m-\ell+1}B_i\mathcal B_{\ell-i}(\lambda^b)\frac{z^m}{m!}\\
&=\frac{1}{\lambda}\sum_{m=0}^\infty\sum_{\ell=0}^{m+1}\sum_{i=0}^\ell\binom{m+1}{\ell}\binom{\ell}{i}\frac{a^{m-l+i+1}b^{m-i+1}}{(m-\ell+2)(m+1)}B_i\mathcal B_{\ell-i}(\lambda^b)\frac{z^m}{m!}\,.
\end{align*}
Comparing the coefficients on both sides of (\ref{eq:gf}), we get the following expression. 
\begin{theorem}  
For $\lambda\ne 0$ with $\lambda^a=1$ and $\lambda^b\ne 1$, and a nonnegative integer $m$, 
\begin{align*}
S_m^{(\lambda)}(a,b)&=\sum_{\ell=0}^{m+1}\sum_{i=0}^\ell\binom{m+1}{\ell}\binom{\ell}{i}\frac{a^{m-l+i+1}b^{m-i+1}}{(m-\ell+2)(m+1)\lambda}B_i\mathcal B_{\ell-i}(\lambda^b)\\ 
&\quad -\frac{\mathcal B_{m+1}(\lambda)}{(m+1)\lambda}\,. 
\end{align*}
\label{th:m1}
\end{theorem}

\noindent 
{\it Remark}  
When $\lambda=-1$ in Theorem \ref{th:m} or Theorem \ref{th:m1}, formulas for Sylvester sums (5.11)--(5.14) in \cite{wang08} are obtained. For, when $a$ is odd, $\mathcal B_n((-1)^a)=-n E_{n-1}(0)/2$ ($n\ge 0$), where $E_n(x)$ are Euler polynomials defined by 
$$
\frac{2e^{x z}}{e^z+1}=\sum_{n=0}^\infty E_n(x)\frac{z^n}{n!}\quad(|z|<\pi)\,.  
$$   
\medskip

In particular, when $\lambda=-1$ and $m=1,2$ in Theorem \ref{th:m1}, we have the following formulas. The first relation is not included in the formula in Theorem \ref{th:1}.  

\begin{Cor}
When $a$ is even and $b$ is odd, 
\begin{align*}
S_1^{(-1)}(a,b)&=\frac{b(a b-a-b)+1}{4}\,,\\
S_2^{(-1)}(a,b)&=\frac{a b(b-1)(2 a b-a-3 b)}{12}\,. 
\end{align*}
\label{cor:m1}
\end{Cor} 

For example, for $a=4$ and $b=11$, we get 
\begin{align*} 
&S_1^{(-1)}(4,11)\\
&=(-1)^0\cdot 1+(-1)^1\cdot 2+(-1)^2\cdot 3+(-1)^4\cdot 5+(-1)^5\cdot 6+(-1)^6\cdot 7\\
&\quad +(-1)^8\cdot 9+(-1)^{9}\cdot 10+(-1)^{12}\cdot 13+(-1)^{13}\cdot 14+(-1)^{16}\cdot 17\\
&\quad +(-1)^{17}\cdot 18+(-1)^{20}\cdot 21+(-1)^{24}\cdot 25+(-1)^{28}\cdot 29\\
&=80\,. 
\end{align*}
From Corollary \ref{cor:m1}, we also get 
$$
S_1^{(-1)}(4,11)=\frac{11(4\cdot 11-4-11)+1}{4}=80\,. 
$$  
Similarly, $S_2^{(-1)}(4,11)=1870$.

\section{Acknowledgement} 

The authors thank the anonymous referees for useful comments and suggestions.

\end{document}